\documentclass[12pt]{amsart}

\def\checkbox{\leavevmode\vbox to 9pt{\hrule \vss
	\hbox to 9pt{\vrule height 9pt \hfil\vrule height 9pt}\vss
	\hrule}\ }

\newcommand{\Q}{\mathbb Q}
\newcommand{\R}{\mathbb R}
\newcommand{\Z}{\mathbb Z}
\newcommand{\N}{\mathbb N}

\renewcommand{\epsilon}{\varepsilon}
\renewcommand{\phi}{\varphi}

\newtheorem{Theorem}{Theorem}[section]

\newtheorem{Definition}{Definition}[section]

\begin{document}
\thispagestyle{plain}
\pagestyle{plain}

\author[]{Alexandr Borisov}
\title[]{ Convex lattice polytopes and cones with few lattice points inside, from a birational geometry viewpoint}
\address{Department of Mathematics, Washington University, St. Louis, MO $63130$, USA}
\email{borisov@math.wustl.edu}
\maketitle


\section{Introduction}

It is pretty well-known that toric Fano varieties of dimension $k$ with terminal singularities correspond to convex lattice polytopes $P \subset \R^k$ of positive finite volume, such that $P \bigcap \Z^k$ consists of the point $0$ and vertices of $P$ (cf., e.g. \cite{Bor-Bor}, \cite{Reid}). Likewise, $\Q-$factorial terminal toric singularities essentially correspond to lattice simplexes with no lattice points inside or on the boundary (except the vertices). There have been a lot work, especially in the last 20 years or so on classification of these objects. The main goal of this paper is to bring together these and related results, that are currently scattered in the literature. We also want to emphasize the deep similarity between the classification of toric Fano varieties and classification of $\Q-$factorial toric singularities.

This paper does not contain any new results. It does contain a sketch of an alternative proof of the qualitative version of the theorem of Hensley (cf. \cite{Hensley}). The paper is organized as follows. In section 2 we discuss some known results about toric Fano varieties, i.e. convex lattice polytopes. In section 3 we discuss some results about $\Q-$factorial toric singularities, i.e. simplicial rational cones. In section 4 we explain some similarities between the above topics, and give a short geometric proof of Hensley's theorem. We also point out some open questions.

We should note that our interest in this subject is motivated by classification questions arising in the Minimal Model Program (cf. \cite{CKM}, \cite{FlipsAbundance}, \cite{Reid}). We assume that the reader is familiar with the basic constructions of the theory of toric varieties. The good reference sources for these are the books of W. Fulton (cf. \cite{Fulton}), T. Oda (cf. \cite{Oda}) and the paper of V. Danilov (cf. \cite{Danilov}). For a discussion of related problems from a purely combinatorial point of view, we refer to the survey by Gritzmann and Wills (cf. \cite{GritzmannWills}). Some good discussion and references can also be found in the book of G. Ewald \cite{Ewald}.

{\bf Notations.} By lattice polytopes of dimension $k$ we will mean closed polytopes of finite positive volume in $\R^k$ whose vertices belong to the lattice $\Z^k \subset \R ^k$. To save space, we will sometimes identify the algebraic geometry objects, like toric Fano varieties, with the corresponding combinatorial objects, like convex lattice polytopes. Hopefully, it will not lead to confusion.

{\bf Acknowledgments.} We thank V. Batyrev who first introduced us to this circle of problems. We also thank G. Sankaran and J.-M. Kantor for helpful discussions.


\section{Toric Fano varieties}

As explained in \cite{Bor-Bor}, cf. also \cite{Fulton}, \cite{Reid} the isomorphism classes of toric Fano varieties $X$ of dimension $k$ are in 1-to-1 correspondence with isomorphism classes of convex lattice polytopes of dimension $k$ with a fixed lattice point $0$ inside. Depending on how bad the singularities of $X$ are allowed to be we have the following sets of equivalence classes of toric Fano varieties, and point-containing convex lattice polytopes.

\begin{itemize}
\item Smooth, to be denoted by $S = S(k)$

\item Terminal, to be denoted by $T = T(k)$

\item Canonical, to be denoted by $C = C(k)$

\item Gorenstein, to be denoted by $G = G(k)$

\item For every $\epsilon$ such that $ 0< \epsilon \le 1,$ $\epsilon -$logterminal, $T_{\epsilon} = T_{\epsilon} (k)$

\item For every $\epsilon$ such that $ 0< \epsilon \le 1,$ $\epsilon -$logcanonical, $C_{\epsilon} = C_{\epsilon} (k)$

\item For every $n\in \N,$ those with Gorenstein index $n,$ $G_n = G_n(k)$

\end{itemize}

Here the singularity is called $\epsilon -$logterminal ($\epsilon -$logcanonical) if its total log-discrepancy is greater than (or equal to) $\epsilon $ (cf. \cite{FlipsAbundance}). Please consult \cite{Bor-Bor}, or \cite{Reid}, or section 3 of this paper for the corresponding combinatorial conditions.

Obviously, $T=T_1,$ $C=C_1,$ $G=G_1.$ From the definitions $S\subseteq T\subseteq C\subseteq T_{\epsilon} \subseteq C_{\epsilon}$ for any $\epsilon \le 1.$ Also, $S\subseteq G_1 \subseteq G_n \subseteq C_{1/n}$ for any $n\in \N.$ In dimension two, some of these classes are the same, because all terminal singularities are smooth and all canonical singularities are Gorenstein. This is true for all singularities, not necessarily toric (cf. e.g. \cite{CKM}).

There are two types of results on classification of toric Fano varieties: the general finiteness theorems and the explicit classification theorems.

The most general finiteness result is given by the following theorem.

\begin{Theorem} (A. Borisov -L. Borisov, 1992, \cite{Bor-Bor}) For any $k\in \N,$ $\epsilon >0,$ the set $C_{\epsilon} (k)$ is finite.
\end{Theorem}

A weaker version of this theorem was proven earlier by V. Batyrev.
\begin{Theorem} (V. Batyrev, 1982, \cite{Batyrev1}) For any $k,n \in \N,$ the set $G_n(k)$ is finite.
\end{Theorem}

It should be mentioned that the above theorems are toric cases of the more general boundedness conjectures for Fano varieties (cf. \cite{Fano}). 
The combinatorial statement that corresponds to Theorem 2.1 is the following.

\begin{Theorem} (D. Hensley, 1983, \cite{Hensley}) For any $k\in \N,$ $\epsilon >0,$ there are only finitely many (up to $GL_k(\Z)$ action) convex lattice polytopes $P$ of dimension $k$ such that $(\epsilon \cdot P )\bigcap \Z^k = \{0\}.$
\end{Theorem}

This theorem was first proven by D. Hensley in 1983 (cf. \cite{Hensley}). Hensley also proved a bound on the volume of such polytopes. This bound was improved by Lagarias and Ziegler (cf. \cite{LagariasZiegler}). The proof of the above theorem in \cite{Bor-Bor} is similar but ineffective. In section 4 of this paper we will sketch a simple geometric proof of it (also ineffective). 
Theorem 2.3 also has the following interesting corollary.

\begin{Theorem} (D. Hensley, 1983, \cite{Hensley}) For any any $k, m \in \N$ there are only finitely many convex lattice polytopes of dimension $k$ with exactly $m$ points strictly inside, up to lattice isomorphisms. 
\end{Theorem}

One must note that in the above theorem $m\ge 1,$ otherwise the statement is false. 

The particular classification theorems are mostly concerned with the sets $S$, $T,$ $C$, and $G$ for small values of $k.$ The smooth case was studied the most.
For $k=2$ the classification is very easy, there are only 5 examples: $P^1\times P^1,$ and $P^2$ with up to three blown-up points. (The points that can be blown up are the closed orbits of the torus action on $P^2$).

 For $k=3$ the classification was done independently by V. Batyrev (cf. \cite{Batyrev3}), and K. Watanabe and M. Watanabe (cf.
\cite{WatanabeWatanabe}). It consists of 18 examples. For $k=4$ the situation is more complicated. There was a lot of work on this, beginning with the thesis of V. Batyrev (cf. \cite{Batyrevthesis}). Unfortunately, \cite{Batyrevthesis} contained some mistakes in the case-by-case analysis which resulted in missing cases. It was partially fixed in \cite{Batyrevpreprint}. The recent preprint of H. Sato \cite{Sato} contains 124 polytopes which is most probably the complete list. For $k\ge 5$ there are some general results, due to V. Batyrev, G. Ewald (cf. \cite{Ewald1}), Voskresenskii and Klyachko (cf. \cite{VoskresenskiiKlyachko}), and others. We refer to the preprints of Batyrev (cf. \cite{Batyrevpreprint}) and Sato (cf. \cite{Sato}), and the book of Ewald (cf. \cite{Ewald}) for explanations and further references.

The set $T(k)$ is known for $k=2$ where $T(2)=S(2)$. For $k=3$ all $\Q-$factorial toric Fano varieties with Picard number 1 were classified by A. Borisov and L. Borisov (cf. \cite{terminal}). Such varieties correspond to lattice tetrahedra with one lattice point inside and no points on edges or faces. There are 8 examples. All but one of the corresponding varieties are weighted projective spaces with the following weights.

(1,1,1,1), (1,1,1,2), (1,1,2,3), (1,2,3,5), (1,3,4,5), (2,3,5,7), (3,4,5,7)

Combinatorially, to obtain the corresponding lattice tetrahedra one can take the tetrahedra in $\R^3$ such that $0$ is a linear combination of vertices with the above coefficients. Then for the lattice one should take the lattice in $\R^3$ generated by the vertices of such tetrahedron.  

One more variety is a quotient of $P^4$ by some action of the group $\Z/5\Z$. This corresponds to taking the lattice tetrahedron that corresponds to $P^4$ and enlarging the lattice to some bigger lattice with relative index 5.

This classification relies on a computer, but its essential part is computer-free. The whole set $T(3)$ was also studied by A. Borisov and L. Borisov with extensive use of computer (cf. \cite{terminal}). We found the minimal and maximal such polyhedra, with respect to the natural embedding ordering. There are 13 minimal and 9 maximal ones. The total list was never produced, because the computational complexity of checking the pairwise non-equivalence of polyhedra was too big for the slow computer that we used at that time. It is expected to contain several hundreds of examples.

The set $C(2)$ is pretty easy to determine. It consists of 16 elements. It was determined, among others by V. Batyrev, cf. \cite{Batyrevthesis}. We refer to \cite{KreuzerSkarke1} for the sketch of an easy proof. We should note that the paper of S. Rabinowitz (cf. \cite{Rabinowitz}) on this topic, referred to in \cite{GritzmannWills} unfortunately misses one example. The set $C(3)$ is probably too big for a reasonable classification. The set of all simplexes in there was determined by A. Borisov and L. Borisov, using a computer (cf. \cite{list}). It contains 225 elements.

The Gorenstein toric Fano varieties are important for mathematical physics. As noticed by Batyrev (cf. \cite{Batyrevmirror}) they provide examples for the mirror symmetry conjecture. For this reason they received considerable attention, especially among physicists. As noticed before, the set $G(2)$ is equal to $C(2)$, so it is known (cf. the paragraph above). The set $G(3)$ was studied M. Kreuzer and H. Skarke, using computer. They found 4319 such polytopes. Kreuzer and Skarke went further to obtain some results for $k=4$. The good account of these results can be found at M. Kreuzer's webpage ``http://tph16.tuwien.ac.at/~kreuzer/CY".


\section{$\Q-$factorial toric singularities}

A $\Q-$factorial toric singularity $X$ of dimension $k$ is just a quotient of the affine space $A^k$ by a finite abelian subgroup of $GL_k$. It corresponds to a simplicial rational cone $C$ in $\R^k =(\Z^k)\otimes \R,$ where $\Z^k$ is the lattice of one-dimensional subgroups of the torus. This cone does not contain any non-trivial linear subspaces of $\R^k.$ We will also assume that it has maximal dimension, otherwise the corresponding variety is (non-canonically) isomorphic to a product of some torus and a lower-dimensional toric singularity (cf., e.g. \cite{Fulton}, \cite{Oda}, \cite{Danilov}).
 
Let us denote by $P_1, P_2, ..., P_k$ the closest to zero lattice points on the extremal rays of $C.$ There is exactly one linear function $\phi$ on $\R^k$ such that $\phi (P_i)=1$ for all $i.$ Because $C$ is rational, $\phi$ has rational coefficients.

The Gorenstein index of the singularity $X(C)$ is equal to $(\phi (\Z^k) : \Z),$ the least common multiple of the denominators of coefficients of $\phi$. The minimal log-discrepancy is the smallest value of $\phi$ on the lattice points in $C.$ Actually, there are two versions of minimal log-discrepancy, the total log-discrepancy and the Shokurov log-discrepancy. The first one is defined using the exceptional divisors of all possible birational morphisms $Y\rightarrow X.$ The second one only uses the divisors whose image is the distinguished point of $X.$ In our case we have such distinguished point, the closed orbit of the torus action. The total log-discrepancy is the smallest non-zero value of $\phi$ on the lattice points in the closed cone $C.$ The Shokurov log-discrepancy is the smallest value of $\phi$ on the lattice points in the open cone $C.$

The total log-discrepancy is obviously not bigger than the Shokurov log-discrepancy, and is also at most 1. the Shokurov log-discrepancy is equal to $k$ for smooth points, and is at most $k/2$ otherwise. Both log-discrepancies are positive, which reflects the fact that $\Q-$factorial toric singularities (and any quotient singularities in general) are log-terminal (cf. \cite{CKM}, \cite{Reid}).

The singularity is called $\epsilon-$logterminal ($\epsilon -$logcanonical) if the total log-discrepancy is greater than (or equal to) $\epsilon$. If $\epsilon =1$ we get just the definitions of terminal (canonical) singularities. These classes of singularities are very important for the Minimal Model Program (cf. \cite{CKM}, \cite{FlipsAbundance}). The most general finiteness result in this area is that for any fixed $k$ and $\epsilon $ all $\epsilon-$logterminal ($\epsilon -$logcanonical) singularities form finitely many ``series" (cf. \cite{vasq}). In order to explain what it means, let us first review the known classification results for terminal and canonical singularities for small $k.$

The terminal singularities is the most restrictive class of singularities that has to be allowed in the Minimal Model Program to make it work. As such terminal singularities have been extensively studied. In dimension 2 there it is easy to see that every terminal singularity is smooth. In dimension 3 the analytic classification exists, in the general case due to S. Mori, M. Reid, and others (cf. \cite{Mori}, \cite{ReidYPG}). A part of it is the classification in the toric case. It is the following theorem, often referred to as ``Terminal Lemma".

\begin{Theorem} (D. Morrison - G. Stevens, 1984, \cite{MorStev}) Every three-dimensional terminal toric singularity is isomorphic to a quotient of $A^3$ by a group $\mu_n$ which acts linearly with weights $\frac1n (1,a,n-a)$ for some $n\in \N$ and $a \in \Z/n\Z,$ with $gcd(a,n)=1$. Here $\mu_n$ is the group of $n-$th roots of unity. The notation $\frac1n (1,a,n-a)$ means that $\rho \in \mu_n$ multiplies the first coordinate by $\rho,$ the second coordinate by $\rho ^a,$ and the third coordinate by $\rho^{(n-a)}.$
\end{Theorem}

In fact, in \cite{MorStev} this theorem is stated for the three-dimensional cyclic quotient singularities. However it is easy to see that every isolated $\Q-$factorial toric singularity is a cyclic quotient. Theorem 3.1 can also be stated as follows. Suppose $x=(x_1, x_2, x_3)$ is a generator of the finite cyclic subgroup of the torus $T^3=\R^3/\Z^3$ that corresponds to the singularity. Then the singularity is terminal if and only if (up to permutation of variables) $x_1+x_2\equiv 0$ mod $\Z$, and for all $k\in \Z$ such that $kx\neq 0$ in $\R^3/Z^3$ none of the coordinates of $kx$ is $0$ in $\R/\Z$. In other words, the subgroup belongs to one of the fixed three 2-dimensional subtori of $T^3$ and intersects trivially with some of their smaller subtori. 

The proof of Theorem 3.1 relies on the combinatorial lemma due to G. K. White, for which D. Morrison and G. Stevens also proposed a new proof. To explain this lemma, and why it is relevant, we first need to explain how $\Q-$factorial terminal toric singularities are related to the lattice-free simplexes. Here by a lattice-free simplex (or, in general, a lattice-free polytope) we will mean a simplex (or polytope) whose vertices are in the lattice, and which contains no other lattice points inside or on the boundary.

In one direction, to any $\Q-$factorial terminal toric singularity one can associate a simplex, which is the set of all points in $x\in C$ such that $\phi (x) \le 1.$ The terminality of the singularity is equivalent to the simplex being lattice-free. In the other direction, for any lattice-free simplex with a distinguished vertex one can construct a rational cone by translating the simplex to put this vertex to the origin and generating the simplicial cone using the other vertices. This cone will determine a $\Q-$factorial terminal toric singularity. Thus, the equivalence classes of the $\Q-$factorial terminal toric singularities are in one-to-one correspondence with the equivalence classes of the lattice-free simplexes with a distinguished vertex. With this in mind, the classification of such singularities in dimension three is equivalent to the theorem 3.2 below. To formulate this theorem we first need the following definition.

\begin{Definition} (cf., e.g. \cite{KannanLovasz}) Suppose $P$ is a convex lattice-free polytope (or, in general, any convex body in $R^k$). Then its width  is the minimum of the lengths of its projections to $\R$ using linear functions on $R^k$ with {\bf integer} coefficients.
\end{Definition}

It is clear that if a convex lattice polytope has width 1, it contains no lattice points inside (though it may still have some on the boundary). The following theorem is a kind of a converse statement for the tetrahedra.

\begin{Theorem} (G. K. White, 1964, \cite{White}). Every lattice-free tetrahedron has width 1.
\end{Theorem}

In fact, G. K. White proved a stronger statement. He allowed the tetrahedra to have lattice points on one pair of opposite edges. There is also the following generalization of Theorem 3.2.

\begin{Theorem} Every 3-dimensional lattice-free polytope has width 1.
\end{Theorem}

We were unable to trace the proof of this theorem. V. Danilov attributes it to M. A. Frumkin (cf. \cite{Danilov1}). H. Scarf attributes it to R. Howe (cf. \cite{Scarf} It looks like neither of the proofs has been published. 

We should note that if one allows points on the boundary, the width may be bigger than 1, even in dimension two. On the other hand, the width is always bounded by a function of dimension. This result is in fact very general, it holds for all convex lattice-free bodies in $\R^k,$ not necessarily lattice polytopes. This theorem is quite old, and there is a long history of improvements on the bound. We refer to the paper of Kannan and Lov\'asz \cite{KannanLovasz} for a very good bound and further references. We should also mention that in higher dimensions there exist lattice-free simplexes of arbitrarily large width, by the result of J.-M. Kantor (cf. \cite{Kantor}).

In dimension 4 much less is known.  D. Morrison and G. Stevens classified all abelian quotients of dimension 4 that are both terminal and Gorenstein (cf. \cite{MorStev}). In 1988 S. Mori, D. Morrison, and I. Morrison used computer to study terminal $(\Z/p\Z)-$quotients, for prime $p$ (cf. \cite{MMM}). The singularities they discovered seemed to fall into the finite number of series, similar to the series of three-dimensional terminal singularities above. More precisely, they found 1 three-parameter series, 2 two-parameter series, 29 one-parameter series and several thousands of 0-parameter series (what they called unstable singularities). Here the series of terminal singularities of dimension 3 is considered to have 2 parameters. They conjectured that their list was complete. It was partially proven in 1990 by G. Sankaran (cf. \cite{Sankaran}). He showed that all the ``stable" four-dimensional prime quotient singularities are among those found in \cite{MMM}. Together with our result (cf. \cite{vasq}) this implies that the Mori-Morrison-Morrison list is complete up to possibly finitely many exceptions. Together with the extensive computer evidence of \cite{MMM}, it is quite likely that there are indeed no other such singularities.

For canonical singularities, the classification is known for $k\leq 3$. For $k=2$ it is very easy to see that they are cyclic quotients of the type $\frac{1}{n}(1,n-1)$ for some $n\in \N$. In dimension 3 it was done by M.-N. Ishida and N. Iwashita (cf. \cite{IshidaIwashita}).
In fact, their result is very general. They obtained a complete classification of all 3-dimensional canonical toric singularities, including those that are not $\Q-$factorial. In the $\Q-$factorial case there are two 2-parameter series, one 1-parameter series and two 0-parameter series (exceptional singularities). We refer to \cite{IshidaIwashita}, Theorem 4.1 for the details. We should also note that similar but weaker results were obtained independently by D. R. Morrison, using different methods (cf. \cite{Morrison}).

In general, for any fixed $k$ and $\epsilon$ the $\epsilon-$logterminal (logcanonical) singularities form a finite number of series. The general definition of a series is somewhat complicated, we refer to \cite{vasq} for the details. The implication for the cyclic quotients of prime index is the following. We formulate it for the  $\epsilon-$logterminal case, the same is true verbatim for the $\epsilon-$logcanonical singularities.

\begin{Theorem} (A. Borisov, 1997, \cite{vasq}) For any fixed $k$ and $\epsilon$ there is a finite collection of closed subgroups $\{V_i\}$ of the torus $T^k=\R^k/\Z^k$ and a finite collection of their closed subgroups $V_{i,j} \subset V_i$ such that the $(\Z/p\Z)-$quotient singularity is $\epsilon-$logterminal if and only if a generator of the corresponding subgroup of $T^k$ belongs to $V_i \backslash (\bigcup \limits_{j} V_{i,j})$ for some i. 
(It is clear that the above condition does not depend on the choice of a generator).
\end{Theorem}

This theorem, and a more general theorem of \cite{vasq} for all $\Q-$factorial toric singularities in fact follows from the theorem of J. Lawrence (cf. \cite{Lawrence}). The good approximation to this theorem, which is sufficient for the theorem 3.5 above, is the following.

\begin{Theorem} (J. Lawrence, 1991, \cite{Lawrence}) Suppose $S$ is a closed subset of a finite-dimensional torus $T$ such that $nS \subseteq S$ for all $n\in N.$ Then $S$ is a finite union of the closed subgroups of $T.$
\end{Theorem}

The main theorem of Lawrence is more general. 
\begin{Theorem} Suppose $T$ is as above (or, possibly, $T$ is is not a torus but a a closed subgroup of some torus). Suppose $U$ is an open subset of $T$ (or, more generally, a {\it full} subset of $T$, i.e. for all closed subgroups $L$ of $T$ the intersection of $L$ and $U$ is either empty or contains a relatively open subset of $L$). Consider all closed subgroups of $T$ that don't intersect $U$. Then the number of maximal elements of it, with respect to inclusion, is finite.
\end{Theorem}

The Lawrence's proof of it is very elegant and well-written, which makes the paper \cite{Lawrence} a must-read for anyone seriously interested in the subject. It uses some geometry of numbers. One can also give a geometric proof of its weaker version (Theorem 3.5 above). It is similar to the geometric proof of the Hensley theorem, which will be discussed in the next section. It is however somewhat complicated, and will possibly appear elsewhere.

\section{Some similarities and open questions}
It was noticed in particular by J. Lawrence that the topics of the above two sections have something in common. Namely, the kind of geometry of numbers involved in the proof of the main theorem in \cite{Lawrence} is similar to what was used by Hensley in \cite{Hensley} (and later by Lagarias -Ziegler, cf. \cite{LagariasZiegler}, and Borisov-Borisov, cf. \cite{Bor-Bor}). Another similarity is the following. The classification of terminal (and canonical) weighted projective spaces by Borisov-Borisov (cf. \cite{terminal}) is very similar to the case-by-case analysis of Sankaran (cf. \cite{Sankaran}).

In our opinion, the main driving force behind both results of Hensley and Lawrence is some elementary properties of the dynamics of multiplication by integers on a torus. To explain this, we will sketch a short conceptual proof of the (qualitative) theorem of Hensley. One can also prove in a similar manner a weak theorem of Lawrence (Theorem 3.6). But it is somewhat complicated, so it will possibly appear elsewhere. We should also note that the same ideas were used in \cite{Toricdiscr} to prove Shokurov's conjecture that minimal discrepancies of toric singularities can only accumulate from above.

\begin{Theorem} (D. Hensley, \cite{Hensley}, also \cite{LagariasZiegler}, \cite{Bor-Bor}) For any $k\in \N$ and $\epsilon >0,$ there are only finitely many (up to $GL_k(\Z)$ action) convex lattice polytopes $P$ of dimension $k$ such that $(\epsilon \cdot P )\bigcap \Z^k = \{0\}.$
\end{Theorem}

{\bf Sketch of the proof.} For simplicity, we will present the proof for the  case $\epsilon =1$. However essentially the same proof works in the general case. First of all, by using Minkowski Lemma one can reduce the problem to the case of simplexes, whose vertices generate the whole lattice (cf., e.g. \cite{Bor-Bor}). If we move the simplex to put one of the vertices to zero, we get two lattices. First one, which we will now call $\Z^k$ is generated by the vertices. The second one is the original lattice. It contains $\Z^k$ and the quotient subgroup is generated by the point $O$ in $P$ which corresponds to the zero of the original lattice. This reduces the problem to showing that for every $k$ there are just finitely many points $O$ in the standard open simplex $\Delta \subset T^k$ such that the finite subgroup of $T^k$ generated by $O$ contains no other points from $\Delta$. 

Suppose there are infinitely many such points. By compactness of $\bar{\Delta}$ we can find an infinite sequence $\{O_i\}$ of such points that converge to some point $O^*\in \bar{\Delta}.$

Suppose first that $O^*\in \Delta.$ Then for some natural $n >1 $ $nO^*$ also belongs to $\Delta.$ If $nO^*\neq O^*$ then for $i$ big enough $O_i$ is close to $O^*$ and $nO_i$ is close to $nO^*$, so they are different points in $\Delta,$ contradiction. If $nO^*=O^*$ then take $\epsilon >0$ such that the ball of radius $\epsilon$ $B_{\epsilon}(O^*)$ is contained in $\Delta$. For $i$ big enough $|O_i-O^*| < \epsilon /n.$ Therefore the point
$$nO_i=nO^*+n(O_i-O^*)=O^*+n(O_i-O^*)$$
is in $\Delta,$ and is different from the point $O_i=O^*+(O_i-O^*)$, contradiction.

Finally, if $O^*$ belongs to the boundary of $\Delta,$ choose the face $\Delta '$ which interior it belongs to. Choose $n>1$ so that $nO^*$ also belongs to the interior of $\Delta '$. Then the same argument as above works, because since $O_i$ approach $O^*$ from the inside of $\Delta,$ the $nO_i$ also approach $nO^*$ from the inside of $\Delta.$ One just needs to choose $\epsilon$ small enough so that $B_{\epsilon}(O^*)$ only intersect the faces that contain $O^*.$
This completes the proof.

We would like to mention now several possible directions of research in this area.

1) Try to classify toric Fano varieties of small dimension with relatively mild singularities, using a computer.

2) Try to understand better the smooth and Gorenstein toric Fano varieties. Consult \cite{Sato} and \cite{Batyrevpreprint} for some particular conjectures.

3) Write a computer code that would find explicitly the finite set of series given in Theorem 3.5. In particular, automatize the argument of Sankaran (cf. \cite{Sankaran}).

4) Try to generalize Theorem 3.5. to the non-$\Q-$factorial case. One problem with this is that it is not exactly clear what to mean by a series of singularities in this more general context. One non-trivial result in this direction is due to M.-N. Ishida and N. Iwashita (cf. \cite{IshidaIwashita}). As a related question, what can be said about the set of all lattice-free polytopes of given dimension? 

5) Try to understand better the set of Shokurov log-discrepancies of $\Q-$factorial toric singularities. We proved in \cite{vasq} that these log-discrepancies can only accumulate from above, and only to such log-discrepancies of smaller dimensions. This proved the toric case of a more general conjecture of Shokurov. It also suggests that ``stable" cyclic quotient singularities somehow come from lower-dimensional singularities. This was also noticed by I. Morrison (cf. \cite{MorrisontoSankaran}). It would be interesting to understand exactly how it happens, and maybe get a better conceptual understanding of terminal cyclic quotients of arbitrary dimension.

6) The effective versions of Hensley's theorem (cf. \cite{Hensley}, \cite{LagariasZiegler}) provide bounds for the volumes of the corresponding polytopes that are asymptotically close to the actually existing examples (cf. \cite{PerlesZaksWills}). It would be very interesting to find an effective version of the theorem of Lawrence, and to determine the asymptotics of the number of series, and other parameters involved. We should note here the paper of J.-M. Kantor, who proved the existence of higher-dimensional lattice-free simplexes with arbitrarily large width (cf. \cite{Kantor}). 

There are many other open questions in the area. Some of them can be found in the survey of Gritzmann and Wills (cf. \cite{GritzmannWills}). We would like to stress that many of the problems and methods involved here are quite elementary. On the other hand, it is related to many very advanced areas of modern mathematics and mathematical physics. This makes it a good starting ground for beginning researchers. We hope that this short survey would help bring some more people into this interesting area.


\end{document}